\documentclass[12pt,reqno]{amsart}
\usepackage{graphicx,color}
\usepackage[centertags]{amsmath}
\usepackage{enumerate}
\usepackage{natbib}

\begin{document}

\title[Equi-probable priors]{On the non-arbitrary assignment of equi-probable priors}

\maketitle
\baselineskip=26pt

\begin{center}
{\bf William M. Briggs}  \\ \vskip .05in General Internal Medicine,
Weill Cornell Medical College \\ 525 E. 68th, Box 46,
New York, NY 10021 \\ \textit{email:} mattstat@gmail.com \\
\vskip .1in
\today
\end{center}

\newpage
\baselineskip=26pt

\vskip .1in\noindent \textsc{Summary:}
How to form priors that do not seem artificial or arbitrary is a central question in Bayesian statistics.  The case of forming a prior on the truth of a proposition for which there is no evidence, and the definte evidence that the event can happen in a finite set of ways, is detailed.  The truth of a propostion of this kind is frequently assigned a prior of 0.5 via arguments of ignorance, randomness, the Principle of Indiffernce, the Principal Principal, or by other methods.  These are all shown to be flawed.  The {\it statistical syllogism} introduced by Williams in 1947 is shown to fix the problems that the other arguments have.  An example in the context of model selection is given.

\vskip .1in\noindent \textsc{Key words:} Induction; Logical Probability; Model selection; Principle of Indifference; Principal Principle; Prior formation.
\newpage

\section{Introduction}
There are (at least) two central foundational problems in statistics: how to justify a given probability model, and how to assign prior probabilities to its parameters. The goal of both of these operations is to insure that they are not seen to be arbitrary, or appear to be guided by the subjective whim of the user, and that they logically follow from the explicit evidence that is given or assumed to be known.

As is well known, concepts such as exchangeability, symmetry, and even direct appeals to physics or biology are sometimes given to answer the first question.  Just as well known, however, is that frequently no justification other than habit---or ignorance of any alternative---is used to guide a user to select a particular model.  I do not attempt to examine model correctness here. I take models as given, and instead look at the second question of prior assignment.

That problem is huge, so here I only take a small piece of it: how to form a prior on the truth of an event in two situations: when nothing is known about event other than that it could happen, and when we know that the event can happen in a certain number of finite ways.

First, it is useful to recall, what is often forgotten, that---both deductive and non-deductive---statements of logic are nothing more than the study {\it between} relations, and only between the relations explicitly stated.  How true is one thing given another? is the usual, and the only, question.  The existence and characteristics of the relations themselves is left to other disciplines.  Forgetting this distinction can lead, and has lead, to unnecessary arguments about the nature of logical probability.  My attempt here, is to clear up some of the controversies in the context of prior formation on the truth of an elementary propostion.  An example in model selection will be given.

The meta-logic I assume to evaluate the arguments of this paper is classical \citep{Sch2005}.  So, suppose $p$ is a premise and $q$ a conclusion to the argument from $p$ to $q$.  We may write this in many ways, but one of the clearest is this:
\begin{equation}
\label{arg}
	\frac{\:\:p\:\:\:\:\:\:\:\:\:\:\:\:}{\:\:q\:\:\:\:\:\:\:\:\:\:\:\:}
\end{equation}
which is to be read, ``(the proposition) $p$ (is true) therefore (the proposition) $q$ (is true)." Logical probability makes statements like this:
\begin{equation}
\label{logic}
	a\le \Pr(q|p) \le b.
\end{equation}
I take, from Cox \citeyearpar{Cox1961}, and like those before him \citep{Lap1996, Jef1998, Key2004, Jay2003}, that the values of the constants are $a=0$ and $b=1$, and that if those limits apply to the conclusion $q$ of a given argument, then $q$ is, respectively, certainly false or certainly true.  When those limits are reached, which is rarely, and even never, in the world of statistics, then the logical connective (between $q$ and $p$) is said to be {\it deductive}. If the limits are not reached, then the logical connective is said to be {\it non-deductive}.  Non-deductive arguments may be inductive, or they may be otherwise.  The arguments from $p$ to $q$ are either valid if they are deductive, or invalid if they are not deductive.  Invalid does {\it not} imply unreasonable; neither does deductive imply reasonable.

Here is a simple example of a deductive argument that is not reasonable (in the sense of relevance) adapted from from \citet{Sch2005}: ``If it is raining now, then red is a color. It is raining now.  Therefore, red is a color."  This is a valid argument in classical logic because the implication is always true, as red is certainly a color.  

Inductive arguments---which are argument from contingent premisses which are, or could have been, observed, to a contingent conclusion about something that has not been, and may not be able to be, observed---are, of course, central to probability.  In an earlier paper \citep{Bri2006}, I started with an example of an inductive argument which everybody believes is reasonable.  That was, ($p=$)Because all the many flames observed before have been hot, that ($q=$) {\it this} flame will be hot. Notice that no measure of reasonableness is given, no measure of how true the conclusion $q$ is with respect to its premiss.  We can give such a measure, and that we can do so is explained using the principles of logical probability (which I do not prove here; but see the references below: if you are not yet convinced that probability as logic is the correct interpretation of probability, pretend that it true is for the course of this paper).

The flames argument is inductive. {\it Not} all non-deductive arguments are inductive.  Carnap \citeyearpar{Car1950}, the undisputed point scorer for logical probability in the 20th century, unfortunately had the habit of calling all non-deductive inferences `inductive', which, among other things, lead to a confusion about what logical probability is, and it is this confusion that is in part responsible for logical probability's current refugee status \citep{Sto1973,Sto1983}.  In any case, I do not follow Carnap's terminology here.

For an example of a common, non-inductive (and non-deductive) argument, suppose $M$ is any non-contradictory contingent (not logically necessary) fact or proposition, and $t$ any tautology. The argument:
\begin{equation}
\label{tautology}
	\frac{\:\:t\:\:\:\:\:\:\:\:\:\:\:\:}{\:\:M\:\:\:\:\:\:\:\:\:\:\:\:}
\end{equation}
is not valid (and is read ``$t$, therefore $M$").  Writing out details in this manner makes clear the tacit process of argumentation that is part of any prior probability assignment: all of our evidence is first amassed and then explicitly laid out {\it before} the probability assignment is made.  We are interested in the probability that $M$ is true given $t$.  The most common tautology used in cases like this is $t=$``$M$ can be true or false." The principles of logical probability gives:
\begin{equation}
\label{nond_logic}
	0< \Pr(M|t) < 1.
\end{equation}
And that is the best we can ever do knowing only $t$ that $M$ is contingent (e.g.,  \citet{Key2004}).  This is a point well worth reflecting upon: and which is amplified below. It follows from the logical principle that it is impossible to argue validly to a contingent conclusion given a necessarily true or tautologous premiss.  This result is not dependent on a particular $t$; any tautology $t'$ will do. This is why statements about the probability of $M$ that lack evidence frequently write (\ref{nond_logic}) as ``$0< \Pr(M) < 1$": the missing tautology, since it can be anything, is implied.  This is usually harmless enough, but it can lead to troubles.

Now, the probability statement (4) represents the best that can be said in the face of no evidence, except for the evidence that we know $M$ is contingent: (4) is, or should be, the prior assigned in the face of true ignorance.  Since this prior is not definite, we cannot move towards definitness unless we learn something more about $M$.

Some statisiticians---of the (subjective) Bayesian persuasion---would not like to settle for (\ref{nond_logic}), which is a  vague enough statement about $M$, and would insist that we find some concrete real number $r$ such that $\Pr(M|t)=r$.  To find this number, there is usually an appeal, to the utterers of (\ref{nond_logic}), to announce some subjective opinion they might have about $M$, or even, if it can be believed, about how they would take bets with the Deity (or, for the secular, with Mother Nature) over the truth of $M$. This line was begun by Ramsey, and is summarized in e.g. \citep{Pre2003}. I find this approach wholly unsatisfying.  And so do those who still call themselves frequentists, and who still do so, at least in part, because of their distaste with the wilfull subjectivity of the `Bayesians' and their insistence on using terms like `gambling' and `betting.'

I now examine some of these arguments to find an $r$, and show how they rely on $t$ and on other evidence.

Not all Bayesians would insist that you must say how you'd bet for or against $M$.  Some try to find $r$ by an argument like the following: ``Well, $M$ can be true, or it may be false.  So it must be that $\Pr(M)=\frac{1}{2}$."  No, it musn't.  The first sentence to this argument is just $t$  (the missing, but implied, conclusion before the probability statement, is just $M$) , and nothing has been gained.  The step from the conclusion to the probability statement is arbitrary (as many have felt before; e.g. \citep{Fis1973b,Goo1983}).

The argument can be modified, and moved away from strict ignorance, by inserting some additional evidence, by, say, $e=$``$M$ is equally like to be true or false", which I hope you agree is the same as saying $e=$``$\Pr(M|t)=\frac{1}{2}$."  The argument is then:

	\vspace{.1in}
	\begin{center}
		\begin{tabular}{p{3.5in}cc} \baselineskip=26pt
			$M$ is true or it is false & & \\
			&&\\
			$\Pr(M|t)=\frac{1}{2}$ & & \\
			&&\\\cline{1-1}
			&&\\
			$\Pr(M|t)=\frac{1}{2}$ .& & \addtocounter{equation}{1}(\theequation)\\
		\end{tabular}
	\end{center}
	\vspace{.1in}

This is a curiously dogmatic argument; nevertheless, it {\it is} a valid one; however, the (major) premiss is the same as the conclusion, which isn't wrong, but it is begging the question.  This is usually and loosely called a fallacy, but the conclusion {\it does} follow from assuming the premisses are true, therefore the argument {\it is} valid: it is just of no use. (A helpful way to read this argument is to say ``$p$ is true, therefore $p$ is true."  Attaching the tautology $t$, or any other tautology or necessary truth, changes nothing; it is then ``$p\& t$ is true, therefore $p\& t$ is true.")

There is still the matter of assigning a probability statement to the conclusion of (5), which is:
\begin{equation}
\label{weird}
	\Pr\left(``\Pr(M|t)=\frac{1}{2})"\:|\:e,t\right)=1,
\end{equation}
a statement which is cruicial to understand: it just says that the conclusion deductively follows from the premisses (the reasoning to (6) is equivalent to that which led Egon Spengler to say, ``There's definitely a very slim chance we'll survive.").

The argument (5) is usually recognized for what it is, and instead, in their search for $r$, people will more likely say ``Well, $M$ can be true, or it may be false, {\it and I have no reason to think that it is false or that it is true. Both are equally likely.}  So it must be that $\Pr(M)=\frac{1}{2}$."  This kind of argument is sometimes called the ``Principle of Indifference," advanced by \citep{Key2004} and criticized in e.g. \cite{HowUrb1993}.  It is the ``indifference" or ``no reason" clause that is start of troubles.

So the primary purpose of this paper is to show what that `no reason' premiss does and does not mean, and how arguments for and against models are influenced by it, and by other evidence.

\section{No Reason}
In argument (5), the phrase ``Both [possibilities for $M$] are equally likely" is no more than a restatement of the premiss ``$\Pr(M|t)=\frac{1}{2}$," which makes whole argument as useless as it was when the premiss was explicitly stated in numerical form.  However, it is evidently itself a conclusion from the premiss, ``I have no reason to think that $M$ is false or that it is true."  Now, this argument, in its many forms, has lead a happy life.  It, or a version of it, shows up in discussion of priors frequently, and also, of course, in discussions about model selection, e.g. \citep{BerSmi2000}.  But it is an argument that should not have had the attention it did.  For we can rewrite it like this :

	\vspace{.1in}
	\begin{center}
		\begin{tabular}{p{3.5in}cc} \baselineskip=26pt
			I do not know---I am {\it ignorant}; I have complete ignorance---whether $M$ is true or false, but it can only be true or false. & & \\
			&&\\\cline{1-1}
			&&\\
			$M$ & & \addtocounter{equation}{1}(\theequation)\\
		\end{tabular}
	\end{center}
	\vspace{.1in}

\noindent The conclusion to (7) is usually assigned probability $\Pr(M|t)=\frac{1}{2}$.  This argument, I hope you can see, is not valid: the conclusion certainly does not follow from the premiss, and the probability statement is arbitrary.  Here's why.  This argument {\it is} valid:

	\vspace{.1in}
	\begin{center}
		\begin{tabular}{p{3.5in}cc} \baselineskip=26pt
			$M$ is true or it is false & & \\
			&&\\\cline{1-1}
			&&\\
			I do not know---I am {\it ignorant}; I have complete ignorance---whether $M$ is true or false, but it can only be true or false. & & \addtocounter{equation}{1}(\theequation)\\
		\end{tabular}
	\end{center}
	\vspace{.1in}

\noindent It should now be obvious that the conclusion is nothing more than a restatement of the initial tautology!  To be explicit: saying you do not know, or are ignorant, about $M$ is making a statment that is equivalent to $t$.  So, despite our repeated insistence of ignorance, we are back to (\ref{nond_logic}), which is to say, right where we started.  It should, therefore, be---but it is not---astounding that people have instead come to the probability statement ``$\Pr(M|t)=\frac{1}{2}$" for the conclusion of (7).  I think I know why this is so.

Up to this point, I have been very careful not to give an example for $M$, some concrete, real-world thing upon which to fix the idea in your mind. This was on purpose.  Because it is difficult, if not nearly impossible, especially if you are a working statistician, to avoid adding hidden premisses to (3) or (7) once you have such an example in mind, and then to criticize the conclusion that (4) is valid.  To emphasize: (4) is the correct statement to make given that the {\it only} evidence for $M$ is $t$ and that $M$ is contingent.

To validly arrive at an $r$, new evidence about $M$ must be added.  These additional premisses have to be of a certain concrete character themselves.  They cannot be anything like ``$M$ can be true or false" or any other restatement of $t$.  They cannot contain the probability statement of the conclusion, as in ``$M$ is equally likely true or false."  Nor can they measure some form of `ignorance,' because that is nothing different than ``$M$ can be true or false."  The best---in the sense of being the most precise---probability statement that can be made given these arguments is (\ref{nond_logic}).  So if we are to find an $r$ what can these premisses be?

Before I tell you, let me first fill in the blank about $M$, and give you a real example.  When I do, unless you are a highly unusual person, you will almost certainly {\it instantly} think, ``Of course the probability of $M$ is a $\frac{1}{2}$!  What is the problem!"

Let $M$ represent the fact that I see a head when next I flip this coin.

Are you with the majority who insist that the probability of $M$ must be $\frac{1}{2}$?  I ask you why this must be so.  But before you answer notice that the `coin flip' $M$ is {\it entirely} different from any other $M$ where all you know is that it is contingent.  For example, if instead of a coin flip, suppose $M$ represented the fact that you open a box and examine some object inside and note whether you can see, without touching anything, an `H'.  Now all you know is that $M$ can be true or false.  Based {\it solely} on the information you have, you do not know any other possibilities.  You do {\it not} know that an `H' or some other letter or object might appear.  You do not know, even, whether a snake may jump out of the box.  If you imply that because the question asked something about an `H' that the result must be `H' or some other letter, probably a `T', then you are {\it adding} evidence that you were {\it not} given.

Back to the coin flip.  Why is the probability of $M$ $\frac{1}{2}$?  Symmetry, perhaps?  As in, ``It can fall head or tail and there is no reason to prefer---I am indifferent---to head over tail"? But isn't that the same as ignorance, that is, the same as the tautology?  It is.  Because substitute `be true' for `fall head' and `be false' for `fall tail' and you are right back at the tautology.  Or symmetry as in, ``Heads and tails are equally likely because I have no reason to think otherwise"?  Again, the ``no reason to think otherwise" is the ignorance argument, and the ``Heads and tails are equally likely" or ``indifference" is begging the question.

The anticlimatic answer is the statistical syllogism, as defined by Williams \citeyearpar{Wil1947} in this example:

	\vspace{.1in}
	\begin{center}
		\begin{tabular}{p{3.5in}cc} \baselineskip=26pt
			Just 1 out of 2 of the possible sides are Heads & & \\
			&&\\
			$M$ is an side & & \\
			&&\\\cline{1-1}
			&&\\
			$M$ is a head & & \addtocounter{equation}{1}(\theequation)\\
		\end{tabular}
	\end{center}
	\vspace{.1in}

\noindent This inductive argument is, of course, invalid.  But we can now justify saying $\Pr(M|e)=\frac{1}{2}$, where $e$ is the two premisses.  Adding arguments to $e$ about symmetry, or `fair' coins, or ignorance does not change the probability of the conclusion, because these arguments are all equivalent to adding $t$ or ``$P(M|t)=1/2$" to the list of premisses.  A `fair' coin, after all, carries with it the assumption that the probability {\it is} 1/2: which is begging the question. 

Another example. Suppose there are 10 men in a room and just 9 of these 10 are Schmenges.  $M$ is a man in the room.  The conclusion ``$M$ is a Schmenge" by the statistical syllogism has probability $\frac{9}{10}$.  Note that you rarely hear the term ``fairness"  applied to situtations like this (is there such a thing as a `fair' room full of Schmenges?).  This may be because when there are more than two possibilities for $M$, people are naturally more suspicious that something other than equi-probability holds.

There is no proof of the correctness of the statistical syllogism: it is taken to be axiomatic.  In fact, it may be considered in this sense to be the primary axiom of assigning logical probability.  I have yet to find anybody who disagress with its truth.  Equally compelling, there is no argument against the statistical syllogism (as described fully in \citep{Sto1973, Sto1983}).  It is true that it gives the same results as the traditional arguments of `ignorance', `fairness', or symmetry give, but it does not carry the same baggage.  The other arguments, while they contain the necessary and sufficient information that $M$ is contingent and has two or ten states, carry hidden assumptions, information that is not explicit and can cause consternation and disagreement, because not everybody would necessarily put the same value on these hidden assumptions.  There is no hidden information to the statistical syllogism.  Except maybe something having to do with ``randomness."

\section{Whiter randomness?}

An objection to the statistical syllogism might have something to do with ``randomness", and how it is invoked to select, or to ``sample", say, men from a room.  This may seem fair line of inquiry because of practical interest in the conclusion $M$.  I may want to bet, say, on the chance that the man I grab is a Schmenge, or there might be other reasons why I want to accurately assess the probability of $M$.  But arguments about randomness are, just as are arguments about ignorance, irrelevant.

If you were to grab a man out of the room randomly: how can you be sure that the probability that he is a Schmenge is $\frac{9}{10}$?  Suppose you were to ``sample" the men by opening the door and grabbing the nearest man and noting whether or not he is a Schmenge.  Or perhaps that doesn't sound ``random" enough to you.  Instead, you order the men inside to polka madly, to run about and bounce off the walls and to not stop; then you reach in a grab one.  This sampling procedure becomes an additional premise, so that we have: 

	\vspace{.1in}
	\begin{center}
		\begin{tabular}{p{3.5in}cc} \baselineskip=26pt
			($e_1$) Just 9 out of the 10 men are Schmenges & & \\
			&&\\
			($e_2$) $M$ is a man in the room & & \\
			&&\\
			($e_3$) Men are arranged in the room randomly & & \\
			&&\\\cline{1-1}
			&&\\
			The man $M$ that I grab will be a Schmenge & & \addtocounter{equation}{1}(\theequation)\\
		\end{tabular}
	\end{center}
	\vspace{.1in}

Here, I take ``randomly" to mean ``I have no idea---I am ignorant---of how the men are arranged".  First, suppose {\it all} we know is that there are men in a room, but {\it nothing} else.  That is, our only evidence is $e_3$, which is just another way of saying, ``There are men in the room, and I have no idea who they are or how they are arranged."  Tacit in this is the idea that there may be some Schmenges in the room, which, of course, means that there may not be any.  That is, $e_3$ is equivalent to, ``$M$ may be true or it may be false".  This is our old friend, the tautology $t$, which we have already seen adds nothing to the argument, or to the probability of the conclusion.

Now, if you {\it did} know something about the way the men were arranged in the room, then this is evidence you {\it must} include in the list of premisses, where it would quite legitimately and naturally change the probability of $M$. But just saying your evidence is ``random", or your experiment was ``randomly" sampled, adds nothing.

This should not be surprising, as Bayesians have long known that randomness is not a concept that is needed in experiments such as patient assignment in clinical trials, e.g. \citet{HowUrb1993}.

It should also not be necessary to say that we do not need to assume anything about infinite ``trials" of men in rooms to arrive at the probability of $M$.  Some ``objective" Bayesians try this kind of argument in an attempt justify their priors by invoking something called the {\it Principal Principle}, which states 
\begin{quotation}
that if the objective, physical probability of a random event (in the sense of its limiting relative-frequency in an infinite sequence of trials) were known to be $r$ and if no other relevant information were available, then the appropriate subjective degree of belief that the event will occur on any particular trial would also be $r$: \citep[p. 240]{HowUrb1993}.
\end{quotation}
Ignoring the fact that we can never know what happens after an infinite amount of time, or that cannot imagine a infinite number of rooms filled with Schmenges, but pretending that we can, the Principal Principle says ``$\Pr\left(M|\Pr(M)=r\right) =r$" (it adds the premiss ``$\Pr(M)=r$" which is taken to be the `objective' or physical probability of $M$), but which we can now see is just begging the question.

\section{An example and a brief note on arguing correctly and completely}

Suppose you are considering $M_1$ and $M_2$ as the only competing models for some situation.  Then, using the statistical syllogism and the logical probability assignments it implies as above, $\Pr(M_1\vee M_2|e_s)= \Pr(M_1|e_s)+\Pr(M_2|e_s)=1$ and $\Pr(M_1|e_s)=\Pr(M_2|e_s)=\frac{1}{2}$, where $e_s$ represents the statistical syllogism evidence.  This is the justification for starting with equal probability in model selection. After $x$ is observed, then it is easy (in principle) to calculate $\Pr(M_1|x,e_s)$ and $\Pr(M_2|x,e_s)$.

It is no surprise that this is the same point reached by appealing to the Principle of Indifference (or even the Principle of Maximum Entropy for a finite number of model choices; \citet{Jay2003}).  The statistical syllogism gives the same answers as the Principle of Indifference, but not by the same route and, again, without the hidden assumptions.  The built-in question-begging of that principle is gone, and there is no appeal to subjectivity, which many find so distasteful.

Arguments against equi-probable priors in discrete problems often center, as they do in, for example, Laplace's Rule of Succession for the probability of seeing the sun rise tomorrow (read \citet[chap. 18]{Jay2003} for a fascinating look at this oft-cited topic), on evidence external to that problem.  That is, certain evidence $e$ is given for the truth of $M$, and a probability is then logically assigned to it.  But the critic has in mind evidence $e'$, which may contain $e$ but also has arguments different than $e$, which would naturally lead to a different probability assignment.  This is overwhelmingly true for Laplace's example.

The naive information Laplace started with was just that the event, the sunrise, could happen or not (the tautology), that just one of the two possible outcomes was a sunrise, and that tomorrow was an outcome; therefore the sun would rise, to which he assigned a prior probability of 1/2.  He then went on to modify this probability using Bayes's theorem, and was subsequently ridiculed by just about every probabilist since.

Why?  Because these probabilists had different evidence in mind: such as their knowledge of astronomy and the physics of planets rotating about the sun, that the earth had been on its journey for longer than Lapace assumed, and on and on.  None of these different-evidence based criticisms against Laplace's (tacit) use of the statistical syllogism to assign a prior are relevant.  Laplace's statment of ``$P(\mbox{sunrise}|e_s)=1/2$" is logically correct.

Along this same line, it is often heard that one must select priors, either on models or parameters, {\it before} seeing the data, lest the data somehow modify your pure `prior' thoughts.  This view is false, at least in the strict logical sense, because whether you apply the statistical syllogism before or after seeing your data it is irrelevant to the probability you assign.  It is based {\it only}, in the case of model selection for example, on the argument $M_1 \vee M_2$ is an outcome etc. The probability assignment ``$\Pr(M|e_s)=1/2$" is true no matter when in time it was made.

Here is another example that has certainly happened to every statistician in one form or another.  Suppose you are with some civilian who knows you are a statisician, and the both of you witness some odd event (say, a bus passes by from which emerges an ice cream cone, which is deposited on your hat).  Your friend will feel an almost moral compulsion to ask, ``What's the probability of {\it that} happening!"  Some strict Bayesians (and others) will argue that the probability is 1, since the event certainly happened.  But your friend clearly implied a different question, one along the lines of ``What's the probability of that happening in circumstances that are similar, and given you didn't expect to see it?"   Well, this probability may or may not be calculable, but its probability after the fact is still the same it was before the fact as long as the ``fact" is not one of the premisses of the argument that lead to the probability assignment.

\section{Conclusion}

Logical probability is a much neglected subject in the statistical community.  The only book in many years to appear on the subject is the (semi-polemical) Jaynes \citeyearpar{Jay2003}.  The Bayesian revolution from the later part of the 20th century, remarkable in many ways, eschewed logical probability and fixed on the idea that probabilities are subjective.

I believe that it is this focus on subjectivity which made statisticians comfortable with words like ``ignorance", ``fair" (though that term pre-dates the revolution), ``no reason", and especially ``gamble", ``indifferent",``betting" and so on.  These terms themselves {\it feel} or are directly subjective; they are words to put your beliefs behind.  And once you have brought in belief, you make it difficult to discover the hidden assumptions behind your belief.  This small paper only attempted to cast light on a few of these hidden assumptions in the simplest possible situations.  It is certainly not a complete answer to the question of how to assign prior probabilities in an objective way.

But the statistical syllogism can clearly be applied to all situations, such as assigning priors on probability model parameters, when those parameters can take a finite number of values or states.  The class of probability models which contain such parameters may or may not be very large, but it is at least not empty, though it of course does not contain the most frequently used probability models, such as those, say, from the exponential family.  But I make no attempt in this paper to justify, or modify, the use of the statistical syllogism in the case where the number of outcomes is countably or uncountably infinte, as in the case of parameters in models like the normal distribution.

\newpage
\bibliographystyle{apalike}
\bibliography{logic}

\end{document}